\documentclass[a4paper, 11pt]{article}

\usepackage[latin1]{inputenc}
\usepackage[english]{babel}
\usepackage[nottoc]{tocbibind}
\usepackage{emptypage}

\usepackage{amsmath}
\usepackage{amssymb}
\usepackage{amsfonts}
\usepackage{mathrsfs}

\usepackage[final]{showkeys}

\usepackage{color}
\usepackage{graphicx}

\usepackage[all]{xy}

\usepackage{cite}
\usepackage{epigraph}

\usepackage{dsfont}

\usepackage{geometry}
\newcommand{\C}{\mathbb{C}}

\newcommand{\supp}{\mathrm{supp}\,}

\newcommand{\R}{\mathbb{R}}

\newcommand{\Ci}{\mathcal{C}}

\newcommand{\D}{\mathscr{D}}

\newcommand{\met}{\mathsf{Met}}
\newcommand{\mor}{\mathrm{Mor}}

\newcommand{\de}{\partial}

\newtheorem{Teo}{Theorem}
\newtheorem{Rem}{Remark}
\newtheorem{Prp}[Teo]{Proposition}
\newtheorem{Cor}[Teo]{Corollary}
\newtheorem{Lmm}[Teo]{Lemma}

\title{Homology by  metric currents}
\author{Samuele Mongodi\footnote{Supported by the Italian project FIRB-IDEAS ''Analysis and Beyond''.}}

\begin{document}

\maketitle

\section{Introduction}

Metric currents are, in a certain sense, a metric analogous of flat currents, therefore are related to the geometry of the space and of their support. In this short note, we try to give some evidence for the previous statement, by showing that the homology which can be defined by means of metric normal currents coincides, on nice enough metric spaces, with the usual singular homology.

The ``differential forms'' used to define metric currents are, in some sense, projections of the metric space to the euclidean space and the metric currents are particular elements of the appropriate dual; this bears a similarity with the usual definition of cohomology, by means of functionals on maps from the euclidean simplexes to the space.

Our main result is contained in Theorem \ref{teo_finale} and Corollary \ref{cor_finale}, asserting the equivalence of the ``metric'' homology and the usual one in locally Lipschitz contractible CW-complexes; our strategy of proof is quite straightforward: we use the geometric property of metric currents (expressed by Proposition \ref{prp_supp_corr}) and a generalization of the cone construction (introduced in \cite{ambrosio1} and already employed in \cite{mongodi1} for a different problem) to verify that the axioms of Eilenberg and Steenrod for homology hold.

Examples of locally Lipschitz contractible spaces are normed spaces, CAT($k$) spaces and Alexandrov spaces (see \cite{yamaguchi1} and \cite{mistuyama}); Alexandrov spaces and CAT($k$) spaces are metric spaces with some geometric conditions which are in some sense resemblant of curvature bounds. We refer to \cite{burburiv} for precise definitions, properties and examples.

Similar results have been hinted every now and then in the recent literature on metric currents, without proof (for instance, in \cite{wegner1}), or have been presented in more elaborate and refined settings, with a resulting increase of the difficulty and technicality of the proofs (see \cite{riedwegschaeppi}, \cite{mitsuishi}).

Given the great interest in the geometry of metric spaces, in the last decades, a detailed and simple proof of such a basic theorem could be helpful to spread further light on the phenomena related to the (Lipschitz) geometry of metric spaces, introducing a concrete possibility of using metric currents as a tool to explore and understand them.

\bigskip

\emph{Aknowledgement: }I would like to thank Simone Borghesi for a fruitful discussion on the uniqueness theorem of Eilenberg and Steenrod.

\section{The setting}

For the definition of metric currents, we refer to \cite{ambrosio1} and \cite{lang1}. In what follows, we use the local theory of metric currents, developed by U. Lang in \cite{lang1}; in this setting, we only consider locally compact metric spaces. In particular, the completeness assumption, needed in \cite{ambrosio1}, can be dropped.

\medskip

Let $X$ be a locally compact metric space. We will denote by $\D^k$ the space of metric $k$-differential forms with locally Lipschitz coefficients and by $N_k$ the space of normal currents with compact support, i.e. currents with finite mass whose boundary has finite mass; we will also denote by $d$ the boundary operator. We note that, on a complete metric space, the space of normal currents with compact support following Lang's definition coincides with the space of normal currents with compact support following Ambrosio and Kirchheim's definition (see the fourth section of \cite{lang1} for a more detailed comparison between the two theories).

We recall the following result about metric currents and their supports (cfr. \cite[Proposition 3.3]{lang1}).

\begin{Prp}\label{prp_supp_corr}Let $T\in\D_m(X)$ and $A\subset X$ be a locally compact set containing $\supp(T)$. Then there is a unique $T_A\in\D_m(A)$ such that
$$T_A(f,\pi)=T(\tilde{f},\tilde{\pi})$$
whenever $(f,\pi)\in\D^m(A), (\tilde{f},\tilde{\pi})\in\D^m(X)$ and $\tilde{f}\vert_A=f$, $\tilde{\pi}_i\vert_A=\pi_i$. Moreover, $\supp(T)=\supp(T_A)$.\end{Prp}

We proceed now to define the homology complex. Given metric space $X$, we can consider the chain complex
$$\ldots\longrightarrow N_k(X)\stackrel{d}{\longrightarrow} N_{k-1}(X)\longrightarrow\ldots\longrightarrow N_1(X)\stackrel{d}{\longrightarrow}N_0(X)\longrightarrow 0$$
where $N_k(X)$ is the space of normal metric currents with compact support, and the associated homology
$$H_k(X)=\frac{\mathrm{Ker}\{d:N_k(X)\to N_{k-1}(X)\}}{\mathrm{Im}\{d:N_{k+1}(X)\to N_k(X)\}}$$
Obviously, if $f:X\to Y$ is a Lipschitz map, we obtain the pushforward operator $f_\sharp:N_k(X)\to N_k(Y)$ for every $k$ and, since  $f_\sharp$ and $d$ commute, we have an induced operator
$$H(f):H_k(X)\to H_k(Y)$$
such that $H(\mathrm{Id})=\mathrm{Id}$ and $H(f\circ g)=H(f)\circ H(g)$. In other words, $H$ is a covariant functor from the category of metric spaces with Lipschitz functions to the category of abelian groups. In what follows we will write $f_*$ instead of $H(f)$.

Moreover, if $A$ is closed subset of $X$,we define $N_k(X,A)$ setting
$$
N_k(X,A)=N_k(X)/N_k(A).
$$

\section{The axioms of Eilenberg and Steenrod}

Since $d:N_k(X)\to N_{k-1}(X)$ sends $N_k(A)$ in $N_{k-1}(A)$ we can consider the relative homology groups $H_k(X,A)$ and we have the long exact sequence of the pair, in the same way of singular homology
$$\ldots H_k(A)\to H_k(X)\to H_k(X,A)\stackrel{d'}{\to}H_{k-1}(A)\to H_{k-1}(X)\to\ldots $$
where $d'$ is an homomorphism of degree $-1$.

\begin{Prp} \label{prp_mayer_vietoris}
Let $\{U,V\}$ be an open covering of $X$, let $i_U, i_V$ be the inclusions of $U\cap V$ in $U$ and $V$ respectively and let $j_U, j_V$ be the inclusions of $U$ and $V$ respectively in $X$. Then the short sequence of chain complexes
$$0\to N_*(\bar U\cap \bar V)\stackrel{(i_U)_*\oplus(i_V)_*}{\longrightarrow}N_*(\bar U)\oplus N_*(\bar V)\stackrel{(j_U)_*-(j_V)_*}{\longrightarrow}N_*(X)\to0$$
is exact.
\end{Prp}
\noindent{\bf Proof: }
Given $T\in N_k(U\cap V)$ with $(i_U)_\sharp(T)=0$, for every form $(f,\pi)\in\D^k(\bar U)$ we have
$$T(f\vert_{U\cap V},\pi\vert_{U\cap V})=0$$
so $T(g,\eta)=0$ for every $(g,\eta)\in\D^k(\bar{U}\cap\bar{V})$, that is $T=0$.

Moreover, if $(j_U)_\sharp(T)=(j_V)_\sharp(S)$, with $T\in N_k(\bar U)$ and $S\in N_k(\bar V)$, then $\supp ((j_U)_\sharp(T))=\supp((j_V)_\sharp(S))\subseteq \bar U\cap \bar V$; this means that $T=(i_U)_\sharp R$ and $S=(i_V)_\sharp R$ with $R\in N_k(\bar U\cap \bar V)$.

Finally, given $T\in N_k(X)$, we can consider a partition of unity subordinated to the covering $\{U, V\}$, $\{\phi_U,\phi_V\}$. The current $T\llcorner\phi_U$ has support contained in $U$, therefore, by Proposition \ref{prp_supp_corr}, there is $S_1\in N_k(\bar U)$ such that $T_\llcorner\phi_U=(j_U)_\sharp S_1$; similarly, there is $S_2\in N_k(\bar V)$ such that $-T\llcorner\phi_V=(j_V)_\sharp S_2$. So, we have that $T=(j_U)_\sharp S_1-(j_V)_\sharp S_2$ and the exactness of the sequence follows. $\Box$

\medskip

By employing the usual techniques of homological algebra and Proposition \ref{prp_mayer_vietoris}, we can now prove the exactness of the Mayer-Vietoris sequence for the homology of normal currents (see e.g. \cite[Theorem 7.19]{switzer1}).

\begin{Prp}\label{prp_excision}
Given a closed subset $A$ of $X$ and an open set $U$ such that $\bar{U}$ is contained in the interior of $A$, we have that the inclusion map $(X\setminus U, A\setminus U)\to (X,A)$ induces an isomorphism in homology.
\end{Prp}
\noindent{\bf Proof: }
The inclusion $(X\setminus U, A\setminus U)\subseteq (X,A)$ induces a morphism $j_\sharp:N_*(X\setminus U, A\setminus U)\to N(X,A)$; with the same technique of restriction used in the proof of the previous Proposition, $j_\sharp$ is showed to induce a isomorphism in homology. $\Box$

\medskip

\begin{Prp}\label{prp_homot_inv}
Two locally Lipschitz-homotopic locally Lipschitz maps 
$$f\sim g:X\to Y$$
induce the same homomorphism in homology
\end{Prp}
\noindent{\bf Proof: }
Let $H:X\times[0,1]\to Y$ be the locally Lipschitz homotopy between $f$ and $g$ and let us define the operator
$$K:N_k(X)\to N_{k+1}(Y)$$
by the following formula
$$K(T)(f,\pi_1,\ldots, \pi_{k+1})=\sum_{i=1}^{k+1}(-1)^{i+1}\int_0^1T\left(f\circ H \cdot \frac{\de \pi_i\circ H}{\de t},\ldots, \hat{\pi_i},\ldots\right)$$
Arguing like in \cite[Proposition 10.2]{ambrosio1} (see also \cite{mongodi1} for another generalization of the same Proposition), we see that if $T\in N_k(X)$, $K(T)$ is also in $N_{k+1}(Y)$ and the following holds
$$d(K(T))=-K((dT))+g_\sharp T-f_\sharp T.$$
Consequently, if $dT=0$, we see that $g_\sharp T-f_\sharp T$ is in the image of $d:N_{k+1}(Y)\to N_k(Y)$, that is $f_*=g_*$ as applications between $H_*(X)$ and $H_*(Y)$. $\Box$

\medskip

\begin{Prp}\label{prp_dim_ax}
If $X$ is a metric space with only one point, we have
$$H_*(X)=\left\{\begin{array}{ll}\mathbb{K}&\textrm{ if }*=0\\
0&\textrm{otherwise}\end{array}\right.$$
where $\mathbb{K}$ is either $\R$ or $\C$.
\end{Prp}
\noindent{\bf Proof: }
The thesis is obvious, as $M_0(X)=\{\alpha \delta_x\ \vert\ \alpha\in\mathbb{K}\}\cong \mathbb{K}$ and $M_j(X)=\{0\}$ for $j>0$. $\Box$

\medskip

\begin{Teo}\label{teo_axioms}The functor $H_*$ defines a homology theory.\end{Teo}
\noindent{\bf Proof: } We refer to \cite{switzer1, eilsteen1} for a complete discussion of the Eilenberg-Steenrod axioms and their variations; we content ourselves with noting that $H_*$ is obviously a functor and that, by the results we proved in the previous pages, it fulfills the homotopy axiom, the excision axiom, the dimension axiom and  the exactness of the long sequence of the pair.

Moreover, it is obvious that $N_*$ is additive on the disjoint union of spaces and that $d$ preserves such structure, hence also $H_*$ is additive. \hfill $\Box$ 

\section{Uniqueness of homology}

Let $\met_2$ be the category whose objects are pairs $(X,A)$ where $X$ is a locally compact metric space and $A$ is a closed subset in $X$ and whose morphisms are  continuous functions. Similarly, let $\met_{2L}$ be the category whose objects are the same of $\met_2$ and whose morphisms are locally Lipschitz maps; obviously, if $(X,A)$ is an object in $\met_2$ (or $\met_{2L}$), also $(X,\emptyset)$ and $(A,\emptyset)$ are, and all the inclusions between them are morphisms of the category.

Let us denote by $\met_2'$ (resp. $\met_{2L}'$) the category whose objects are the same of $\met_2$ and whose morphisms are equivalence classes of morphisms of $\met_2$ (resp. $\met_{2L}$) with respect to the relation of being homotopic (resp. being homotopic by a locally Lipschitz homotopy).

We will denote by $\mor_{\Ci}(O_1, O_2)$ the class of morphisms from $O_1$ to $O_2$ in the category $\Ci$.

\begin{Lmm}\label{lmm_CWlip} If $(X,A)$ and $(Y,B)$ are objects of $\met_2'$, then they are also objects of $\met_{2L}'$. Moreover, if $(X,A)$ and $(Y,B)$ are locally Lipschitz contractible CW-pairs, then 
$$\mor_{\met_2'}((X,A), (Y,B))=\mor_{\met_{2L}'}((X,A),(Y,B))$$
\end{Lmm}
\noindent{\bf Proof: } The first assertion is obvious, by the definition of the categories $\met_2'$ and $\met_{2L}'$. If $(X,A)$ and $(Y,B)$ are CW-pairs, $X$ and $Y$  are locally finite CW-complexes, as they are metrizable; given a continuous map of pairs $f:(X,A)\to(Y,B)$, by the cellular approximation theorem, we can find a cellular map $f':(X,A)\to(Y,B)$.

Let $e$ be a cell of $X$; its characteristic map is homotopic to a Lipschitz map, by the hypothesis of local Lipschitz contractibility. Let us divide the boundary of $e$ in small simplices such that their images through the characteristic maps are contained in Lipschitz contractible balls; using the procedure described in the beginning of the proof of \cite[Proposition 1.3]{yamaguchi1}, we costruct a Lipschitz map which is homotopic to the characteristic map.

So, given the map $j:e\to X$, if $f'\circ j:e\to Y$ is Lipschitz on the boundary of $e$, we can construct (in the same way as before) a Lipschitz map $g_j:e\to Y$ which coincides with $f'\circ j$ on the boundary of $e$. As $j$ is a homeomorphism on the interior of $e$, this permits us to inductively construct a locally Lipschitz map $g:(X,A)\to  (Y,B)$ which is homotopic to $f'$, hence to $f$. Moreover, if $f'$ was already locally Lipschitz on a sub complex $C\subset X$, we can take $g$ to agree with $f$ on $C$.

\medskip

On the other side, if $f_1,\  f_2:(X,A)\to (Y,B)$ are locally Lipschitz maps, with a continuous homotopy $K$ between them (hence representing the same morphism in $\met_2'$), we can apply the preceding reasoning to the map
$$K:(X,A)\times (I,{0,1})\to (Y,B)\;.$$
As $f_1,\ f_2$ are already locally Lipschitz, we can construct a map $K'$, homotopic to $K$, locally Lipschitz, and such that $K'(x, 0)=f_1(x)$ and $K'(x,1)=f_2(x)$ for every $x\in X$; therefore, $f_1$ and $f_2$ represent the same morphisms also in $\met_{2L}'$. \hfill $\Box$

\bigskip
Lemma \ref{lmm_CWlip} tells us that there is a natural isomorphism between $\met_2'$ and $\met_{2L}'$.

We denote by $H^L_*$ the Lipschitz simplicial homology functor (with real coefficients) on $\met_2'$, i.e. the homology (with real coefficients) obtained by Lipschitz simplicial chains. By \cite[Proposition 1.3]{yamaguchi1}, on a locally Lipschitz contractible space, $H^L_*$ coincides with the usual simplicial homology.%; a CW-complex which is also a complete metric space is therefore locally finite, hence it is locally Lipschitz contractible.

\begin{Teo}\label{teo_finale} There exists a natural transformation from the Lipschitz simplicial homology (with real coefficients) $H^L_*$ to $H_*$, which induces an isomorphism between them on every CW-pair in $\met_2$.
\end{Teo}
\noindent{\bf Proof: } We define a transformation from the chain complex of Lipschitz chains to normal currents with compact support. 

Let $\sigma:\Delta_k\to X$ be a Lipschitz $k$-simplex, let $[e_k]$ be the usual integration current on $e_k$; $\sigma(\Delta_k)$ is compact in $X$, hence the $k$-current $T_\sigma=\sigma_{\sharp}[e_k]$ has compact support. By the properties of metric currents, $T_\sigma$ is a normal current; moreover, $[e_n]$ is also a classical current, hence (see for instance \cite{federer1}) $d[e_n]=[\de e_n]$, where $\de e_n$ is the boundary of $e_n$ in the sense of simplicial chains.

Therefore, $dT_{\sigma}=T_{d\sigma}$. This means that $T$ induces a natural transformation $T_*$ between homologies.

Such a transformation induces an isomorphism between $H_*(X)$ and $H_*^L(X)$ when $X=\{x\}$. Therefore, by the classical results on homology (see \cite[Chapter 7]{switzer1}), the two homology theories are isomorphic, when restricted to CW-pairs in $\met_{2L}'$, hence in $\met_{2}'$, hence in $\met_2$. \hfill $\Box$

\medskip

The following corollary is an obvious consequence of the previous theorem and of the result on locally Lipschitz contractible spaces mentioned above.

\begin{Cor}\label{cor_finale}The homology $H_*$ defined in this note coincides with the usual simplicial homology on locally Lipschitz contractible CW-pairs in $\met_2$.\end{Cor}

\medskip

\begin{Rem} The same argument can be repeated verbatim for the integral homology, employing the space $I_k(X)$ of integral currents and noticing that, by the Boundary Rectifiability Theorem \cite[Theorem 8.6]{ambrosio1}, $dI_k(X)\subseteq I_{k-1}(X)$.\end{Rem}

\medskip

Instead of employing local metric currents, it is possible to adjust the definition of metric current given in \cite{ambrosio1} to consider also non complete metric spaces: if one adds the requirement for the mass of a current to be a tight measure, then one can dispense with completeness, at least as long as the result we employed in this note are concerned. This allows us to obtain the same results for non-complete metric space, as it is done in \cite{mitsuishi}; with this strategy, however, the proofs become much more involved.

\bibliography{bibsing}{}
\bibliographystyle{siam}

\end{document}